\documentclass[12pt,a4paper]{article}
\usepackage{amsfonts}

\usepackage{amssymb}

\usepackage{amsmath}

\newtheorem{theorem}{Theorem}[section]

\newtheorem{lemma}[theorem]{Lemma}

\date{}

\begin{document}
\title{ Gr\"{o}bner-Shirshov Bases for Commutative Algebras
with Multiple Operators and Free Commutative    Rota-Baxter
Algebras\footnote{Supported by  the Young Project on Natural Science
Fund of Zhanjiang  Normal University (No.QL0902).} }

\author{
{\small Jianjun Qiu}\\
{\small \ Mathematics and Computational  Science School }\\
{\small \ Zhanjiang  Normal University}\\
{\small \ Zhanjiang  524048, China}\\
{\small \ jianjunqiu@126.com}
 }

\maketitle \noindent\textbf{Abstract:} In this paper, the
Composition-Diamond lemma for commutative  algebras with multiple
operators is established.  As  applications, the
Gr\"{o}bner-Shirshov bases and linear bases of  free commutative
Rota-Baxter algebra, free commutative $\lambda$-differential algebra
and free commutative $\lambda$-differential Rota-Baxter algebra are
given, respectively. Consequently,  these three free algebras are
constructed directly by commutative $\Omega$-words.



\section{Introduction}

Let $K$ be a unitary commutative ring and $\lambda\in K$. A
Rota-Baxter algebra of weight $\lambda$ is a  $K$-algebra $R$ with a
linear operator $P:R\rightarrow R$ satisfying the Rota-Baxter
relation:
$$
 P(u)P(v)=P(uP(v))+P(P(u)v)+ \lambda P(uv),\forall u, v\in
R.
$$
The Rota-Baxter operator  first occurred in the paper of G. Baxter
\cite{Bax} to solve an analytic problem, and the algebraic study of
this operator was started by G.-C. Rota \cite{ro}.

There have been some constructions of free (commutative) Rota-Baxter
algebras. In this aspect, G.-C. Rota \cite{ro} and P. Cartier
\cite{ca} gave the explicit constructions of free commutative
Rota-Baxter algebras of weight $\lambda=1$, which they  called
shuffle Baxter and standard Baxter algebras, respectively. Recently
L. Guo and W. Keigher \cite{gk, gk1} constructed the free
commutative Rota-Baxter algebras (with unit or without  unit) for
any $\lambda \in K$ using the mixable shuffle product. These are
called the mixable shuffle product algebras, which generalize the
classical construction of shuffle product algebras. K. Ebrahimi-Fard
and L. Guo \cite{EG08a} further constructed the free associative
Rota-Baxter algebras by using the Rota-Baxter words. In
\cite{EG08b}, K. Ebrahimi-Fard and  L. Guo   use rooted trees and
forests to give explicit construction of free noncommutative
Rota-Baxter algebras on modules and sets.

 A differential algebra of weigh $\lambda$, also called
$\lambda$-differential algebra is a $K$-algebra $R$ with a linear
operator $D:R\rightarrow R $ such that
$$
D(uv)=D(u)v+uD(v)+\lambda D(u)D(v), \forall u, v\in R.
$$
Such an operator $D$ is called a $\lambda$-differential operator. E.
Kolchin \cite{ko} considered  the differential algebra and
 constructed
 free   differential algebra of weight  $\lambda = 0$.  L. Guo and W. Keigher
 \cite{lw} dealt with a generalization  of this algebra and using
 the same way to construct the free differential algebra of
  weight $\lambda$ in both commutative and associative case.

Similar to the relation of integral and differential operators, L.
Guo and W. Keigher \cite{lw} introduced the notion of
$\lambda$-differential Rota-Baxter algebra which is a  $K$-algebra
$R$ with  a $\lambda$-differential operator $D$ and a Rota-Baxter
operator $P$ such that $DP=Id_{R}$. In the same paper \cite{lw},
they construct the free $\lambda$-differential Rota-Baxter
(commutative and associative).

The Gr\"{o}bner-Shirshov bases theory for Lie algebras was
introduced by A. I. Shirshov \cite{b09,S3}. Shirshov \cite{S3}
defined the composition of two Lie polynomials, and proved the
Composition lemma for the Lie algebras. L. A. Bokut \cite{b76}
specialized the approach of Shirshov to associative algebras, see
also Bergman \cite{b}. For commutative polynomials, this lemma is
known as the Buchberger's Theorem in \cite{bu65, bu70}.

  The  multi-operators algebras
($\Omega$-algebras) were  introduced  by A. G. Kurosh \cite{ku} and
the Gr\"{o}bner-Shirshov bases for $\Omega$-algebras were given in
the paper by V. Drensky and R. Holtkamp \cite{dr}.
Composition-Diamond lemma for associative algebras with multiple
linear operators (associative $\Omega$-algebras) is  established in
a recent paper by L. A. Bokut, Y. Chen and J. Qiu \cite{b08}. Also,
the Gr\"{o}bner-Shirshov bases for Rota-Baxter algebras is
established by L. A. Bokut, Y. Chen and X. Deng \cite{b10} and the
Composition-Diamond lemma for $\lambda$-differential associative
algebras with multiple operators is constructed by J. Qiu and Y.
Chen \cite{qc}.

In this paper, we deal with  commutative  algebras with multiple
linear operators. We construct free commutative algebras with
multiple linear operators and   establish the Composition-Diamond
lemma for such algebras. As applications, we obtain
Gr\"{o}bner-Shirshov bases of free commutative Rota-Baxter algebra,
commutative $\lambda$-differential algebra and commutative
 $\lambda$-differential Rota-Baxter algebra, respectively. Then, by
 using the Composition-Diamond lemma,
linear bases of these three free algebras are obtained respectively.

The author would like to express his deepest gratitude to Professors
L. A. Bokut and Yuqun Chen for their  kind guidance, useful
discussions and enthusiastic encouragements.

\section{ Free  commutative algebras with multiple operators}
In this section, we construct free  commutative algebras with
multiple linear operators.

 Let $K$ be a unitary commutative ring.
A commutative algebra with multiple operators is a commutative
$K$-algebra $R$ with a set
 $\Omega$ of multi-linear  operators.

Let $X$ be a set, $CS(X)$ the free commutative semigroup on $X$ and
$$
\Omega=\bigcup_{t=1}^{\infty}\Omega_{t}
 $$
where $\Omega_{t}$ is the set of $t$-ary operators.

Define
$$
\Gamma_0= X,\  \mathcal {Q}_{0}=CS(\Gamma_0)
$$
and
$$
\Gamma_{1}= X\cup \Omega(\mathcal {Q}_{0}),\ \mathcal
{Q}_{1}=CS(\Gamma_1)
$$
where
$$\Omega(\mathcal
{Q}_{0})=\bigcup\limits_{t=1}^{\infty}\{\theta(u_1,u_2,\dots,u_t)|\theta\in
\Omega_t, u_i\in \mathcal {Q}_{0}, \ i=1,2,\dots,t\}.
$$
For $n>1$, define
$$
 \Gamma_{n}= X\cup \Omega(\mathcal {Q}_{n-1}),\
\mathcal {Q}_{n}=CS(\Gamma_{n})
$$
where
$$\Omega(\mathcal
{Q}_{n-1})=\bigcup\limits_{t=1}^{\infty}\{\theta(u_1,u_2,\dots,u_t)|\theta\in
\Omega_t, u_i\in \mathcal {Q}_{n-1}, \ i=1,2,\dots,t\}.
$$

Let
$$
\mathcal {Q}(X)=\bigcup_{n\geq0}\mathcal {Q}_{n}.
$$
Then it is easy to see that $\mathcal {Q}(X)$ is a commutative
semigroup such that $ \Omega(\mathcal {Q}(X))\subseteq \mathcal
{Q}(X). $

Let   $K[ X; \Omega]$ be the commutative $K$-algebra spanned by
$\mathcal {Q}(X)$.  Extend linearly each  $\sigma\in\Omega_t$,
$t\geq 1$,
$$\sigma:\mathcal {Q}(X)^t\rightarrow \mathcal {Q}(X), \
(x_1,x_2,\dots,x_n)\mapsto\sigma(x_1,x_2,\dots,x_n)
$$
to
$$
K[ X; \Omega]^t\rightarrow K[ X; \Omega].
$$

 Then, it is easy to see that $K[ X; \Omega]$ is
a free commutative  algebra with multiple operators $\Omega$ on set
$X$.

\section{ Composition-Diamond lemma for   commutative  algebras with
multiple operators}

In this section, we introduce the notions of Gr\"{o}bner-Shirshov
bases  for the commutative algebras with multiple operators and
establish the Composition-Diamond lemma for such algebras.

The element in $\mathcal {Q}(X)$ and $K[ X; \Omega]$ are called
commutative $\Omega$-word  and commutative
  $\Omega$-polynomial, respectively.  For any  $u\in \mathcal {Q}(X)$, $u$ has a unique
expression
$$
u=u_1u_2\cdots u_n
$$
where each  $u_i\in X\cup \Omega(\mathcal {Q}(X))$. If this is the
case, we define $bre(u)=n$.  Let $u\in \mathcal {Q}(X)$. Then
$$
dep(u)=\mbox{min}\{n|u\in\mathcal {Q}_{n} \}
$$
is called the depth of  $u$.

Let $\star\notin X$.  By a commutative $\star$-$\Omega$-word we mean
any expression in $\mathcal {Q}(X\cup \{\star\})$ with only one
occurrence of $\star$.  We define the set of all commutative
$\star$-$\Omega$-words on $X$  by $\mathcal {Q}^\star (X)$. Let $u$
be a commutative $\star$-$\Omega$-word and $s\in K[ X;\Omega]$.
Then we call
$$
u|_{s}=u|_{\star\   \mapsto s} \ \mbox{(replace  $\star$
 by $s$ in $u$)}
$$
is a commutative  $s$-$\Omega$-word.  For example, if
$$
u=x_1\theta_3(x_2, \theta_2(\star x_4 ,x_5), x_6)\in\mathcal
{Q}^\star(X),
$$
then
$$
u|_{s}=u|_{\star\ \mapsto s}=x_1\theta_3(x_2, \theta_2(sx_4 ,x_5),
x_6)
$$
is a commutative  $s$-$\Omega$-word.

Similarly, we can define  commutative  $(\star_1,
\star_2)$-$\Omega$-words as expressions in  $\mathcal {Q}(X\cup
\{\star_1, \star_2\})$  with only one occurrence of $\star_1$ and
only one occurrence of $\star_2$. Let us denote by $\mathcal
{Q}^{\star_1, \star_2} (X)$ the set of all  commutative $(\star_1,
\star_2)$-$\Omega$-words. Let $u\in \mathcal{Q}^{\star_1, \star_2}
(X)$, $s_1, s_2\in K[X; \Omega]$. Then we call

$$
u|_{s_1,\ s_2}= u|_{\star_1\mapsto s_1,\star_2\mapsto s_2},
$$
a commutative  $s_1$-$s_2$-$\Omega$-word.\\

Now, we assume that $ \mathcal {Q}(X)$ is equipped with a monomial
order $>$. This means that $>$ is a well order on $\mathcal {Q}(X)$
such that for any $w, v\in \mathcal {Q}(X), u\in \mathcal
{Q}^\star(X)$,
$$
w>v\Rightarrow u|_{w}>u|_{v}.
$$

Note that such an order on $\mathcal {Q}(X)$ exists, for example,
the order (\ref{o1}) in the next section.

For every commutative $\Omega$-polynomial $f\in K[ X;\Omega]$, let
$\bar{f}$ be  the leading term  of $f$. If the coefficient of
$\bar{f}$ is $1$, then we call $f$  monic.

Let $f, g$ be two monic $\Omega$-polynomials.  Then there are two
kinds of compositions.
\begin{enumerate}
\item[(I)]If there exists a commutative $\Omega$-word $w=a\bar{f}=b\bar{g}$ for some $a,b\in
\mathcal {Q}(X)$ such that $bre(w)< bre(\bar{f})+bre(\bar{g})$, then
we call $(f,g)_{w}=af-bg$ the intersection composition of $f$ and
$g$ with respect to $w$.
\item[(II)] If there exists a commutative $\Omega$-word $w=\bar{f}=u|_{\bar{g}}$ for some
$u \in \mathcal {Q}^\star(X)$, then we call $(f,g)_{w}=f-u|_{g}$ the
including composition of $f$ and $g$ with respect to $w$.
\end{enumerate}

In the above definition, $w$ is called the ambiguity of the
composition. Clearly,
$$
(f,g)_w\in Id(f,g) \ \ \ \mbox{ and }  \ \ \ \overline{(f,g)_w}< w
$$
where $Id(f,g)$ is the ideal of $K[ X;\Omega]$ generated by $f$ and
$g$.

Let $f, g $ be commutative $\Omega$-polynomials and $g$ monic with
$\bar{f}=u|_{\bar{g}}$ for some  $u\in \mathcal {Q}^\star(X)$. Then
the transformation
$$
f\rightarrow f-\alpha u|_{g}
$$
is called  elimination of the leading commutative $\Omega$-word
(ELW) of $f$ by $g$, where $\alpha$ is the coefficient of the
leading commutative $\Omega$-word of $f$.

Let $S$ be a set of monic commutative $\Omega$-polynomials. Then the
composition $(f,g)_w$ is called trivial modulo $(S,w)$, if
$$
(f,g)_w=\sum\alpha_iu_i|_{s_i}
$$
where each $\alpha_i\in K$,  $u_i\in \mathcal {Q}^\star(X)$, $s_i\in
S$ and $u_i|_{\overline{s_i}}< w$. If this is the case, we write
$$
(f,g)_w\equiv 0 \ \ mod (S,w).
$$
In general, for any two commutative $\Omega$-polynomials $p$ and
$q$, $ p\equiv q \ \ mod (S,w)$ means $ p-q=\sum\alpha_iu_i|_{s_i}$
where each $\alpha_i\in K$,  $u_i\in \mathcal {Q}^\star(X)$, $s_i\in
S$ and $u_i|_{\overline{s_i}}< w$.

Then $S$ is called a Gr\"{o}bner-Shirshov basis  in  $K[ X;\Omega]$
if any composition $(f,g)_w$ of $f,g\in S$  is trivial modulo
$(S,w)$.

\begin{lemma}\label{l3.3}
Let $S$ be a Gr\"{o}bner-Shirshov  basis  in  $K[ X;\Omega]$
 and $u_1, \ u_2\in \mathcal {Q}^\star(X), \  s_1, s_2\in S$. If
 $w=u_1|_{\overline{s_1}}=u_2|_{\overline{s_2}}$, then
$$
u_1|_{s_1}\equiv u_2|_{s_2} \ mod (S,w).
$$
\end{lemma}
{\bf Proof:} \ There are three cases to consider.

(I)\ \ The commutative $\Omega$-words $\overline{s_1}$ and
$\overline{s_2}$ are disjoint. Then there exits a commutative
$(\star_1, \star_2)$-$\Omega$-words $\Pi$ such that
$$\Pi|_{\overline{s_1},\
\overline{s_2}}=u_1|_{\overline{s_1}}=u_2|_{\overline{s_2}}.
$$
Then
\begin{eqnarray*} u_2|_{ s_2}-u_1|_{
s_1}&=&\Pi|_{\overline{s_1}, \ s_2}-\Pi|_{s_1, \ \overline{s_2}}\\
&=&(-\Pi|_{s_1-\overline{s_1}, \ s_2}+\Pi|_{s_1, \
s_2})+(\Pi|_{s_1,\ s_2-\overline{s_2}}-\Pi|_{s_1, \ s_2})\\
&=&-\Pi|_{s_1-\overline{s_1}, \ s_2}+\Pi|_{s_1,\ s_2-\overline{s_2}}
\end{eqnarray*}
Let $$-\Pi|_{s_1-\overline{s_1}, \
s_2}=\sum\alpha_{2_t}u_{2_t}|_{s_2}\ \ \mbox{and}\ \ \Pi|_{s_1,\
s_2-\overline{s_2}}=\sum\alpha_{1_l}u_{1_l}|_{s_1}.$$
 Since
$\overline{s_1-\overline{s_1}}<\overline{s_1}$ and
$\overline{s_2-\overline{s_2}}<\overline{s_2}$, we have
$$
u_{2_t}|_{\overline{s_2}}, \ u_{1_l}|_{\overline{s_1}}<w.
$$
Therefore
$$
u_{_{2}}|_{ s_2}-u_1|_{
s_1}=\sum\alpha_{2_t}u_{2_t}|_{s_2}+\sum\alpha_{1_l}u_{1_l}|_{s_1}
$$
with each $u_{2_t}|_{\overline{s_2}}, \
u_{1_l}|_{\overline{s_1}}<w$. It follows that

$$
u_1|_{s_1}\equiv u_2|_{s_2} \ mod (S,w).
$$

(II) The commutative $\Omega$-words  $\overline{s_1}$ and
 $\overline{s_2}$ have nonempty intersection but do not
 include each other. For example,
 $$
a\overline{s_1}=b\overline{s_2}
 $$
for some commutative  $\Omega$-words $a,\ b$. Then there exists a
commutative $\star$-$\Omega$-word $\Pi$ such that
$$
\Pi|_{a\overline{s_1}}=u_1|_{\overline{s_{_1}}}=u_2|_{\overline{s_2}}=\Pi|_{b\overline{s_2}}.
$$
Then we have
$$
u_{_2}|_{ s_{_2}}-u_1|_{ s_{_1}}=\Pi|_{bs_{_2}}-\Pi|_{as_{_1}}
=-\Pi|_{as_{_1}-bs_{_2}}.
$$
Since  $S $ is a Gr\"{o}bner-Shirshov basis in  $K[ X;\Omega]$,
 we have
$$
as_1-bs_2=\sum\alpha_jv_j|_{s_j}
$$
where each $\alpha_j\in K, \ v_j\in \mathcal {Q}^\star(X), \ s_j\in
S$ and $v_j|_{\overline{s_j}}<a\overline{s_1}$. Let
$$
\Pi|_{v_j|_{s_j}}=\Pi_{j}|_{s_j}.
$$
Then
$$
u_2|_{ s_2}-u_1|_{ s_1}=\sum\alpha_j\Pi_{j}|_{s_j}
$$
with
$$
\Pi_{j}|_{\overline{s_j}}<w.
$$
It follows that

$$
u_1|_{s_1}\equiv u_2|_{s_2} \ mod (S,w).
$$

(III) One of commutative $\Omega$-words $\overline{s_1}$,
 $\overline{s_2}$ is contained in the other.
For example, let
$$
\overline{s_1}=u|_{\overline{s_2}}
$$
for some commutative  $\star$-$\Omega$-word $u$. Then
$$
w=u_2|_{\overline{s_2}}=u_1|_{u|_{\overline{s_2}}}
$$
and
$$ u_2|_{ s_2}-u_1|_{
s_1}=u_1|_{u|_{s_2}}-u_1|_{ s_1}=-u_1|_{ s_1-u|_{s_2}}.
$$
Similarly to  (II), we can obtain the result. \hfill
$\blacksquare$\\

The following theorem is an analogy  of  Shirshov's composition
lemma for Lie algebras \cite{S3}, which was specialized to
associative algebras by Bokut \cite{b76}.  For commutative algebras,
this lemma is known as the Buchberger's Theorem in \cite{bu65,
bu70}.

\begin{theorem}\label{3.4}{\em(Composition-Diamond lemma)}\ \  Let $S$ be a set of monic
commutative $\Omega$-polynomials in $K[ X;\Omega]$ and  $>$ a
monomial order on $\mathcal {Q}(X)$.  Then the following statement
are equivalent:
 \begin{enumerate}
\item[(I)] $S $ is a Gr\"{o}bner-Shirshov basis in $K[ X;\Omega]$.
\item[(II)] $ f\in Id(S)\Rightarrow \bar{f}=u|_{\overline{s}}$
for some $u \in \mathcal {Q}^\star(X)$ and $s\in S$.
\
\item[(III)] $Irr(S) = \{ w\in \mathcal {Q}(X) |  w \neq
u|_{\overline{s}}
 \mbox{ for  any} \ u \in \mathcal {Q}^\star(X) \ \mbox{and } s\in S\}$
is a $K$-basis of $K[ X;\Omega]/Id(S)=K[ X;\Omega|S]$.
\end{enumerate}
\end{theorem}
{\bf Proof:} (I)  $\Longrightarrow$ (II)\ \ Let  $0\neq f\in Id(S)$.
Then

$$
f=\sum\limits_{i=1}^{n}\alpha_i u_i|_{ s_i}
$$
where each $\alpha_i\in K$, $u_i\in \mathcal {Q}^\star(X)$ and $
s_i\in S $. Let $w_i=u_i|_{\overline{ s_i}}$ and we arrange this
leading commutative $\Omega$-words in non-increasing order by
$$
w_1= w_2=\cdots=w_m >w_{m+1}\geq \cdots\geq w_n.
$$

We prove the result by  induction  on $m$.

If $m=1$, then $\bar{f}=u_1|_{\overline{ s_1}}$.

Now we assume that $m\geq 2$. Then
$$
u_1|_{\overline{s_1}}=w_1=w_2=u_2|_{\overline{s_2}}.
$$
We prove the result by induction on $w_1$. If  $\bar{f}=w_1$, there
is nothing to prove. Clearly,   $w_1>\bar{f}$. Since $S $ is a
Gr\"{o}bner-Shirshov  basis  in $K[ X;\Omega]$, by Lemma \ref{l3.3},
we have
$$
u_2|_{ s_2}-u_1|_{ s_1}=\sum\beta_jv_j|_{s_j}
$$
where $\beta_j\in K, \ s_j\in S, v_j\in \mathcal {Q}^\star(X)$  and
$v_j|_{\overline{s_j}}<w_1$. Therefore, since
$$
 \alpha_1u_1|_{ s_1}+\alpha_2u_2|_{
s_2}=(\alpha_1+\alpha_2)u_1|_{ s_1}+\alpha_2(u_2|_{ s_2}-u_1|_{
s_1}),
$$
we have
$$
f=(\alpha_1+\alpha_2)u_1|_{ s_1}+\sum\alpha_2\beta_jv_j|_{s_j}+
\sum\limits_{i=3}^{n}\alpha_iu_i|_{ s_i}.
$$

If either $m>2$ or $\alpha_1+\alpha_2\neq 0$, then the result
follows from  induction on $m$. If $m=2$ and $\alpha_1+\alpha_2=0$,
then  the result follows from  induction on $w_1$.

(II)$\Longrightarrow$ (III) For any $f\in K[ X;\Omega]$, by ELWs, we
 can obtain
 that $f+Id(S)$ can be expressed by the elements of $Irr(S)$.
 Now suppose $\alpha_1u_1+\alpha_2u_2+\cdots\alpha_nu_n=0$ in
$ K[X;\Omega|S]$ with $u_i\in Irr(S)$, $u_1>u_2>\cdots>u_n$ and
$\alpha_i\neq 0$.  Then, in $ K[ X;\Omega]$,
$$
g=\alpha_1u_1+\alpha_2u_2+\cdots+\alpha_nu_n\in Id(S).
$$
By (II), we have $u_1=\bar{g}\notin Irr(S)$, a contradiction. So
$Irr(S)$ is $K$-linearly independent. This shows that $Irr(S)$ is a
$K$-basis of $ K[ X;\Omega|S]$.

(III)$\Longrightarrow $(II) Let $0\neq f\in Id(S)$. Suppose that
$\bar{f}\in Irr(S)$. Then
$$
f+Id(S)=\alpha (\bar{f}+Id(S))+ \sum \alpha_i(u_i+Id(S))
$$
where $u_i\in Irr(S)$ and $\bar{f}>u_i$. Therefore, $f+Id(S)\neq 0$,
a contradiction. So $\bar{f}=u|_{\bar{s}}$ for some $s\in S$
 and  $u\in \mathcal {Q}^\star(X)$.

 (II)$\Longrightarrow $(I) \ By the definition of the   composition,
we have  $(f,g)_w\in Id(S)$. If $(f,g)_w\neq 0$, then by (II), $
\overline{(f,g)_w}=u_{1}|_{\overline{s_1}}<w$ for some $s_1\in S$
and $u_1\in \mathcal {Q}^\star(X)$. Let
$$
h=(f,g)_w-\alpha_1u_{1}|_{s_1}
$$
where $\alpha_1$ is the coefficient of   $ \overline{(f,g)_w}$. Then
$\bar{h}<w$ and $h\in Id(S)$. By induction on $w$, we can get the
result. \hfill $\blacksquare$

\section{ Gr\"{o}bner-Shirshov bases for  free commutative Rota-Baxter algebras}

In this section, a Gr\"{o}bner-Shirshov basis for free commutative
Rota-Baxter algebra was obtained. By using the Composition-Diamond
lemma (Theorem \ref{3.4}),  a linear basis of such algebra was given
and  the free commutative Rota-Baxter algebra was directly
constructed by commutative $\Omega$-words.

First of all,  we define an order on $\mathcal {Q}(X)$, which will
be used in  this section. Let $X$ and $\Omega$ be well ordered. We
define an order on $\mathcal {Q}(X)=\bigcup_{n\geq0}\mathcal
{Q}_{n}$ by induction on $n$. For any $u,v \in \mathcal
{Q}_0=CS(X)$, we have
$$
u=x_1^{i_1}x_2^{i_2}\cdots x_t^{i_t}\  \mbox{and} \
v=x_1^{j_1}x_2^{j_2}\cdots x_t^{j_t}
$$
where  each $x_i\in X, \ i_k, j_l\geq 0$ and $x_{i}>x_{i+1}$. Then
we define
$$
 u>v\Longleftrightarrow (bre(u), i_1, i_2, \dots,
i_t)>(bre(v), j_1, j_2, \dots, j_t)\ \ \mbox{lexicographically}.
$$
Assume the order on $\mathcal {Q}_{n-1}$ has been defined. Now, we
define an order on $\Gamma_{n}=X\cup \Omega(\mathcal {Q}_{n-1})$.
Let $v_1, v_2 \in \Gamma_{n}$. Then $v_1>v_2$ means one of the
following holds:

(a) $v_{1},v_{2}\in X$ and  $v_{1}>v_{2}$;

 (b) $v_{1}\in \Omega(\mathcal {Q}_{n-1})$ and $v_{2}\in X$;

(c) $v_{1}=\theta_k(v_{1}')$, $v_{2}=\theta_l(v_{2}'), v_{1}'\in
\mathcal {Q}_{n-1}^k, v_{2}'\in \mathcal {Q}_{n-1}^l $ with
$\theta_k>\theta_l$ or $\theta_k=\theta_l$, $v_{1}'>v_{2}'$
lexicographically.

 For any $u,v \in \mathcal {Q}_n$, we have
$$
u=w_1^{i_1}w_2^{i_2}\cdots w_t^{i_t}\  \mbox{and} \
v=w_1^{j_1}w_2^{j_2}\cdots w_t^{j_t}
$$
where each $w_i\in\Gamma_{n}$,  $i_k, j_l\geq 0$ and
$w_{i}>w_{i+1}$. Here, $bre(u)= i_1+i_2+\cdots+i_t$. Define
\begin{equation}\label{o1}
 u>v\Longleftrightarrow (bre(u), i_1, i_2, \dots,
i_t)>(bre(v), j_1, j_2, \dots, j_t)\ \ \mbox{lexicographically}
\end{equation}
Then the Order  (\ref{o1}) is a monomial order   on $\mathcal
{Q}(X)$.

Let $K$ be a commutative ring with unit  and $\lambda\in K$. A
commutative Rota-Baxter algebra of weight $\lambda$ ( see \cite{Bax,
gk, ro}) is a commutative $K$-algebra $R$ with a linear operator $
P:R\rightarrow R$ satisfying the Rota-Baxter relation:
$$
 P(u)P(v) = P( P(u)v) +P(uP(v))+\lambda P( uv), \  \forall u,v \in  R.
$$

It is obvious that any commutative Rota-Baxter algebra is a
commutative algebra with multiple operators  $\Omega$, where
$\Omega=\{P\}$.

In this section, we assume that $\Omega=\{P\}$. Let $\mathcal
{Q}(X)$ be defined as before with  $\Omega=\{P\}$ and  $K[ X;P]$ be
the free commutative algebra with operator $\Omega=\{P\}$ on set
$X$.

\begin{theorem}\label{t4.1}
With the order (\ref{o1})  on $\mathcal {Q}(X)$,
$$
S=\{P(u)P(v) - P( P(u)v) - P(uP(v))-\lambda P(uv) |\  u,v \in
\mathcal {Q}(X)\}
$$
is a Gr\"{o}bner-Shirshov basis in $K[ X;P]$.
\end{theorem}
{\bf Proof:} The ambiguities of  all possible compositions of the
commutative $\Omega$-polynomials in $S$ are:
$$
(i) \ \ \ \  P(u)P(v)P(w) \ \ \ \ \ \ \ \ (ii) \ \ \ \
P(z|_{P(v)P(w)})P(u)\ \  \ \
$$
where $u, v, w\in \mathcal {Q}(X), z\in \mathcal {Q}^\star(X)$. It
is easy to check that all these compositions are trivial. Here,  for
example, we just check $(i)$. For any $ u, v\in\mathcal {Q}(X)$, let
$$
f(u, v)=P(u)P(v) - P( P(u)v) - P(uP(v))-\lambda P(uv).
$$
Then
\begin{eqnarray*}
&&(f(u,v), f(v,w))_{P(u)P(v)P(w)}\\
&&= -P( P(u)v)P(w) - P(uP(v))P(w)-\lambda P(uv)P(w)\\
&& \ \  \  +P(u) P( P(v)w) + P(u)P(vP(w))+\lambda P(u) P(vw)\\
&&\equiv -P(P( P(u)v)w)-P(P(u)v P(w))-\lambda P(P(u)vw)\\
&& \ \  \ -P(P(uP(v))w)-P(uP(v)P(w))-\lambda P(uP(v)w)\\
&& \ \  \ -\lambda P(P(uv)w)-\lambda P(uvP(w))-\lambda^2 P(uvw)\\
&& \ \  \ +P(P(u) P(v)w) + P(uP( P(v)w))+\lambda P(uP(v)w)\\
&& \ \  \ +P(P(u)vP(w))+P(uP(vP(w)))+\lambda P(uvP(w))\\
&& \ \  \ +\lambda P(P(u)vw)+\lambda P(uP(vw))+\lambda^2P(uvw)\\
&&\equiv -P(P( P(u)v)w)-P(P(uP(v))w)-\lambda P(P(uv)w)\\
&& \ \ \ -P(uP(P(v)w))-P(uP(vP(w)))-\lambda P(uP(vw)) \\
&& \ \ \ +P(P(P(u) v)w)+P(P(uP(v))w)+\lambda P(P(uv)w)\\
&& \ \ \ +P(uP(P(v)w))+P(uP(vP(w)))+\lambda P(uP(vw))\\
 &&\equiv 0 \ \ mod(S,P(u)P(v)P(w)). \ \ \ \ \ \ \
\ \ \ \ \ \ \ \ \  \ \ \ \ \ \ \ \  \  \  \ \   \ \ \ \ \ \ \ \ \ \
\ \ \ \blacksquare
\end{eqnarray*}
Define
\begin{eqnarray*}
\Phi_0&=&CS(X),\\
\Phi_{1}&=&\Phi_0\cup P(\Phi_0)\cup \Phi_0P(\Phi_0)  ,\\
\vdots\ \ & &\ \ \ \ \vdots\\
\Phi_{n}&=&\Phi_0\cup P(\Phi_{n-1})\cup\Phi_0P(\Phi_{n-1}), \\
\vdots\ \ & &\ \ \ \ \vdots
\end{eqnarray*}
where
$$
\Phi_0P(\Phi_{n-1})=\{ uP(v)|u\in \Phi_0, v\in \Phi_{n-1}\}.
$$

Let
$$
\Phi(X)=\bigcup_{n\geq 0 }\Phi_n.
$$
The elements in $\Phi(X)$ are called commutative Rota-Baxter words.

By  Theorem \ref{t4.1} and Theorem \ref{3.4}, we have the following
theorems.
\begin{theorem}{\em(\cite{b10})}\label{t4.2}
$Irr(S)= \Phi(X)$  is a basis of  $K[ X;\Omega|S]$.
\end{theorem}

\begin{theorem}
$K[ X;\Omega|S]$ is a free commutative Rota-Baxter algebra  of
weight $\lambda$ on  set $X$ with  a basis $\Phi(X)$.
\end{theorem}

By using  ELWs, we have the following algorithm. In fact, it is an
algorithm to  compute the product of two commutative Rota-Baxter
words in the   free commutative Rota-Baxter  algebra $K[
X;\Omega|S]$.

 \noindent
\textbf{Algorithm} Let $u,v\in  \Phi(X)$. We define $u\diamond v$ by
induction on $n=dep(u)+dep(v)$.

 (1) If $n=0$, then $u,v \in CS(X)$ and $u\diamond
v=uv$.

 (2) If $n\geq 1$, then $u\diamond v=$  \newline
\begin{equation*}
 \left\{
\begin{array}{l@{\quad\quad}l}
uv &  \mbox{if}\  u\in CS(X)\  \mbox{or} \  v\in CS(X)  \\
u_1v_1(P(P(u^{'})\diamond v^{'}) +P(u^{'}\diamond P(v^{'})) +\lambda
P(u^{'}\diamond v^{'})) & \mbox{if} \  u=u_1P(u^{'}),
  v=v_1 P(v^{'})%
\end{array}%
\right.
\end{equation*}

\section{Gr\"{o}bner-Shirshov bases for free commutative $\lambda$-differential
algebras}\label{s3}

In this section, we give  a Gr\"{o}bner-Shirshov basis for a free
commutative $\lambda$-differential algebra and then by   using the
Composition-Diamond lemma (Theorem \ref{3.4}), we obtain a linear
basis of such an algebra, which is the same as the one in \cite{lw}.
Consequently, we construct the free $\lambda$-differential algebra
on set $X$ directly by commutative $\Omega$-words.

 Let $K$ be
a commutative unitary ring and $\lambda\in K$.
 A  commutative $\lambda$-differential  algebra (\cite{lw, ko})
 over $K$  is a  commutative $K$-algebra  $R$
  together with a  linear operator $D:R\rightarrow R$ such that
$$
D(uv)=D(u)v+uD(v)+\lambda D(u)D(v),\ \forall u, v \in R.
$$
  It is obvious that any
commutative $\lambda$-differential
 algebra is a   commutative  algebra with multiple operators $\Omega$ on $X$, where
 $\Omega=\{D\}$.

In this section, we assume that $\Omega=\{D\}$. Let $K[X;D]$ be the
free commutative algebra with operator $\Omega=\{D\}$.   Let $X$ be
well ordered. For any $u \in \mathcal {Q}(X)$, $u$ has a unique
expression
$$
u=u_1u_2\cdots u_n
$$
where $n\geq 1$ and each $u_i\in X\cup D(\mathcal {Q}(X))$.
 Define $ deg(u)$ the number of
$x\in X$  in $u$. For example, if $u=x_1D(D(x_2))D(x_3)$, then
$deg(u)=3$.

We define an order on $\mathcal {Q}(X)=\bigcup_{n\geq0}\mathcal
{Q}_{n}$ by induction on $n$. For any $u,v \in \mathcal
{Q}_0=CS(X)$, we have
$$
u=x_1^{i_1}x_2^{i_2}\cdots x_t^{i_t}\  \mbox{and} \
v=x_1^{j_1}x_2^{j_2}\cdots x_t^{j_t}
$$
where  each  $x_i\in X$, $i_k, j_l\geq 0$ and  $x_{i}>x_{i+1}$. Then
we define
$$
 u>v\Longleftrightarrow ( deg(u), i_1, i_2, \cdots,
i_t)>( deg(v), j_1, j_2, \cdots, j_t)\ \ \mbox{lexicographically}.
$$

 Assume the order on $\mathcal {Q}_{n-1}$ has been defined.
Now, we define an order on $\Gamma_{n}=X\cup \Omega(\mathcal
{Q}_{n-1})$ firstly. Let $v_1, v_2 \in \Gamma_{n}$. Then $v_1>v_2$
means one of the following holds:

(a) $v_{1},v_{2}\in X$ and  $v_{1}>v_{2}$;

 (b) $v_{1}\in D(\mathcal {Q}_{n-1})$ and $v_{2}\in X$;

(c) $v_{1}=D(v_{1}')$,  $v_{2}=D(v_{2}')$ and $v_{1}'>v_{2}'$.

 For any $u,v \in \mathcal {Q}_n$, we have
$$
u=w_1^{i_1}w_2^{i_2}\cdots w_t^{i_t}\  \mbox{and} \
v=w_1^{j_1}w_2^{j_2}\cdots w_t^{j_t}
$$
where each $w_i\in\Gamma_{n}$,  $i_k, j_l\geq 0$ and
$w_{i}>w_{i+1}$. Define
\begin{equation}\label{o2}
 u>v\Longleftrightarrow ( deg(u), i_1, i_2, \cdots,
i_t)>( deg(v), j_1, j_2, \cdots, j_t) \ \ \ \mbox{lexicographically}
\end{equation}

 Then the Order  (\ref{o2}) is a
monomial order   on $\mathcal {Q}(X)$.

\begin{theorem}\label{t5.1}
With the Order  (\ref{o2}) on $\mathcal {Q}(X)$,
$$
S=\{ D(uv)-D(u)v-uD(v)-\lambda D(u)D(v)|\  u,v \in \mathcal {Q}(X)\}
$$
is a Gr\"{o}bner-Shirshov basis in $K[ X;D]$.
\end{theorem}
{\bf Proof:} The ambiguities of  all possible compositions of the
commutative  $\Omega$-polynomials in $S$ are
$$
 D(u|_{_{D(xy)}}v)
$$
where $x, y, v\in \mathcal {Q}(X), u\in \mathcal {Q}^\star(X)$. Let
$$
g(u,v)=D(uv)-D(u)v-uD(v)-\lambda D(u)D(v).
$$

Then
\begin{eqnarray*}
&&(g(u|_{_{D(xy)}}, v),g(x,y))_{D(u|_{_{D(xy)}}v)}\\
&&=-D(u|_{_{D(xy)}})v-u|_{_{D(xy)}}D(v)+\lambda
D(u|_{_{D(xy)}})D(v)\\
&& \ \  \  +D(u|_{_{D(x)y}}v)+D(u|_{_{xD(y)}}v)+\lambda D(u|_{_{D(x)D(y)}}v) \\
&&\equiv -D(u|_{_{D(x)y}})v-D(u|_{_{xD(y)}})v-\lambda D(u|_{_{D(x)D(y)}})v\\
&& \ \  \ -u|_{_{D(x)y}}D(v)-u|_{_{xD(y)}}D(v)-\lambda u|_{_{D(x)D(y)}}D(v)\\
&& \ \  \ -\lambda D(u|_{_{D(x)y}})D(v)-\lambda D(u|_{_{xD(y)}})D(v)-\lambda^2 D(u|_{_{D(x)D(y)}})D(v)\\
&& \ \  \ +D(u|_{_{D(x)y}})v + u|_{_{D(x)y}}D(v)+\lambda D(u|_{_{D(x)y}})D(v)\\
&& \ \  \ +D(u|_{_{xD(y)}})v+u|_{_{xD(y)}}D(v)+\lambda D(u|_{_{xD(y)}})D(v)\\
&& \ \  \ +\lambda D(u|_{_{D(x)D(y)}})v+\lambda  u|_{_{D(x)D(y)}}D(v)+\lambda^2 D(u|_{_{D(x)D(y)}})D(v)\\
&&\equiv 0 \ \ mod(S,D(u|_{_{D(xy)}}v). \ \ \ \ \ \ \ \ \ \ \ \ \ \
\ \  \ \ \ \ \ \ \ \  \  \  \ \   \ \ \ \ \ \ \ \ \ \ \ \ \ \ \ \ \
\ \ \ \ \ \ \ \ \ \ \ \ \ \ \ \  \  \ \ \ \ \ \ \  \ \ \ \ \ \ \
\blacksquare
\end{eqnarray*}

 Let $D^{\omega}(X)=\{D^i(x)|i\geq
0, x\in X\}$, where $D^0(x)=x$  and $CS(D^{\omega}(X))$ the free
 commutative semigroup  generated by $D^{\omega}(X)$.

\begin{theorem}\label{t5.2}{\em(\cite{lw})}
$Irr(S)= CS(D^{\omega}(X))$  is a $K$-basis of $K[X;D|S]$.
\end{theorem}
{\bf Proof:} By Theorem \ref{3.4} and Theorem \ref{t5.1}.

\begin{theorem} {\em(\cite{lw})}
$K[X;D|S]$ is a free commutative $\lambda$-differential  algebra on
set $X$ with  a  basis $CS(D^{\omega}(X))$.
\end{theorem}
{\bf Proof:} By Theorem \ref{t5.2}.\\

By using  ELWs, we have the following algorithm.

\noindent  \textbf{Algorithm} \label{al}(\cite{lw}) Let
$u=u_1u_2\cdots u_n$, where each $u_{k}\in D^{\omega}(X), n>0$.
Define $D(u)$ by induction on $n$.

 (1) If $n=1$, i.e., $u=D^i(x)$ for some $i\geq 0, x\in X $,
 then $D(u)=D^{(i+1)}(x)$.

 (2) If $n\geq 1$, then
$$
D(u)=D(u_1u_2\cdots u_n )=D(u_1)(u_2\cdots u_n )+u_1D(u_2\cdots u_n)
+\lambda D(u_1)D(u_2\cdots u_n).
$$

\section{ Gr\"{o}bner-Shirshov bases for free commutative $\lambda$-differential Rota-Baxter  algebras}

In this section, we give a Gr\"{o}bner-Shirshov basis for a free
commutative $\lambda$-differential Rota-Baxter algebra on a set. By
using the Composition-Diamond lemma for commutative  algebras with
multiple operators  (Theorem \ref{3.4}), we obtain a linear basis of
a  free commutative $\lambda$-differential Rota-Baxter algebra on a
set. Also, we construct the free commutative $\lambda$-differential
Rota-Baxter algebra  on a  set directly by commutative
$\Omega$-words.

Let $K$ be a unitary  commutative ring  and $\lambda\in K$. A
commutative $\lambda$-differential Rota-Baxter algebra (\cite{lw})
is a commutative $K$-algebra $R$ with two linear operators
$P,D:R\rightarrow R$ such that, for any $u,v\in R$,

 (I) (Rota-Baxter
relation) $P(u)P(v)=P(uP(v))+P(P(u)v)+\lambda P(uv),$

(II) ($\lambda$-differential relation) $ D(uv)=D(u)v+uD(v)+\lambda
D(u)D(v),$

(III) $D(P(u))=u$.

It is obvious that any  commutative  $\lambda$-differential
Rota-Baxter  algebra is a commutative algebra with multiple
operators $\Omega$  where
 $\Omega=\{P, D\}$.

In this section, we assume that $\Omega=\{P,D\}$.  Let  $K[ X;\Omega
]$ be  the free commutative algebra with multiple operators $\Omega$
on $X$,  where $\Omega=\{P,D\}$.

Let $X$  be  well ordered and $D>P$. For any $u \in \mathcal
{Q}(X)$, define $ deg_{_{P}}(u)$ the number of $P$ in $u$.

 Define an order on
$\mathcal {Q}(X)=\bigcup_{n\geq0}\mathcal {Q}_{n}$ by induction on
$n$.  For any $u,v \in \mathcal {Q}_0=CS(X)$, we have
$$
u=x_1^{i_1}x_2^{i_2}\cdots x_t^{i_t}\  \mbox{and} \
v=x_1^{j_1}x_2^{j_2}\cdots x_t^{j_t}
$$
where each   $x_i\in X$, $i_k, j_l\geq 0$ and $x_{i}>x_{i+1}$. Then
we define
$$
 u>v\Longleftrightarrow ( deg(u), i_1, i_2, \cdots,
i_t)>( deg(v), j_1, j_2, \cdots, j_t) \ \ \ \
\mbox{lexicographically}.
$$
Assume the order on $\mathcal {Q}_{n-1}$ has been defined. For any
$u,v \in \mathcal {Q}_{n}$,  we have
$$ u=u_1^{k_1}\cdots
u_s^{k_s}\ \mbox{and} \ v=u_1^{l_1}\cdots u_s^{l_s}
$$
where each $u_i \in X\cup \Omega(\mathcal {Q}_{n-1})$, $k_i,l_i \geq
0$ and $u_{i}>u_{i+1}$. Here, $u_{i}>u_{i+1}$ means one of the
following holds:

(a) $u_{i},u_{i+1}\in X$ and  $u_{i}>u_{i+1}$,

 (b) $u_{i}\in D(\mathcal {Q}_{n-1})$ or $u_{i}\in P(\mathcal {Q}_{n-1})$ and $u_{i+1}\in X$,

(c) $u_{i}=\theta_1(u_{i}')$, $u_{i+1}=\theta_2(u_{i+1}')$ where
$\theta_1,\theta_2\in \{D,P\}$ and
$$
 (deg(u_i),deg_{_{P}}(u_i),
\theta_1, u_{i}')> (deg(u_i),deg_{_{P}}(u_i), \theta_2, u_{i+1}')\ \
\ \mbox{lexicographically}
$$
Then we  define $u>v$ if and only if
\begin{equation}\label{o3}
(deg(u),deg_{_{P}}(u), bre(u),k_1,\dots,k_s) >
(deg(v),deg_{_{P}}(v), bre(v),l_1,\dots,l_s)\  \
\mbox{lexicographically}
\end{equation}
Then the Order (\ref{o3}) is a monomial order on $\mathcal {Q}(X)$.

 Let $S$ be a set consisting of  the
following commutative  $\Omega$-polynomials:
\begin{enumerate}
\item[1]
\ $P(u)P(v)-P(uP(v))-P(P(u)v)-\lambda P(uv),\ \  u,v \in \mathcal
{Q}(X)$;
\item[2] $D(uv)-D(u)v-uD(v)-\lambda D(u)D(v), \ \ \  u,v \in \mathcal {Q}(X)$;
 \item[3]$D(P(u))-u, \ \ \ u\in \mathcal {Q}(X).$
\end{enumerate}

\begin{theorem}\label{t6.1}With the Order (\ref{o3}) on
$\mathcal {Q}(X)$, $S$ is a Gr\"{o}bner-Shirshov basis in $K[
X;\Omega]$.
\end{theorem}
\noindent {\bf Proof.}  Denote by $i\wedge j$ the composition of
$\Omega$-polynomials of type $i$ and type $j$. The ambiguities of
all possible compositions of commutative  $\Omega$-polynomials  in
$S$  are only as below. In the following list, $i\wedge j \ \ \ w$
means
 $w$ is the  ambiguity of the composition $i\wedge j$.
 \begin{tabbing}
 $3\wedge3\ \ \ D(P(u|_{_{D(P(v))}}))$,\qquad \qquad \qquad \qquad \qquad \=
 $3\wedge2\ \ \ D(P(u|_{_{D(vw)}}))$\\[0.7ex]
 $3\wedge1\ \ \ D(P(u|_{_{P(v)P(w)}}))$ \qquad \qquad \qquad \qquad \qquad \>
 $2\wedge3\ \ \ D(u|_{_{D(P(v))}}w)$ \\[0.7ex]
 $2\wedge2\ \ \ D(u|_{_{D(vw)}}z)$,\qquad \qquad \qquad \qquad \qquad \>
 $2\wedge1\ \ \ D(u|_{_{P(v)P(w)}}z)$\\[0.7ex]
 $1\wedge3\ \ \ P(u|_{_{D(P(v))}})P(w)$\qquad \qquad \qquad \qquad \qquad \>
 $1\wedge2\ \ \ P(u|_{_{D(vw)}})P(z)$\\[0.7ex]
 $1\wedge1\ \ \ P(z)P(v)P(w)$,\qquad \qquad \qquad \qquad \qquad\>
 $1\wedge1\ \ \   P(v)P(u|_{_{P(w)P(z)}})$
 \end{tabbing}
where $u\in \mathcal {Q}^\star(X), v,w,z\in \mathcal {Q}(X)$. It is
easy to check that all these compositions are trivial.
 Similar to the proofs in Theorem \ref{t4.1} and Theorem \ref{t5.1},
 $1\wedge1$ and $2\wedge2$ are trivial. Others are also easily
checked. Here we just check one  for examples.
\begin{eqnarray*}
2\wedge3&=& -D(u|_{_{D(P(v))}})w-u|_{_{D(P(v))}}D(w)-\lambda
D(u|_{_{D(P(v))}})D(w) +D(u|_{_{v}}w)\\
&\equiv& -D(u|_{_{v}})w-u|_{_{v}}D(w)- \lambda D(u|_{_{v}})D(w)\\
&& +D(u|_{_{v}})w+u|_{_{v}}D(w)+\lambda D(u|_{_{v}})D(w)\\
 &\equiv& 0 \  mod (S,D(u|_{_{D(P(v))}}w)).\ \ \ \ \ \ \ \ \ \ \
\ \ \ \ \ \ \ \ \ \ \ \ \ \ \ \ \ \ \ \ \ \ \ \ \ \  \ \ \ \ \ \  \
\ \ \ \ \ \ \ \blacksquare
\end{eqnarray*}

 Let $D^{\omega}(X)=\{D^i(x)|i\geq
0, x\in X\}$, where $D^0(x)=x$. Define
\begin{eqnarray*}
\Upsilon_0&=&CS(D^{\omega}(X)),\\
\Upsilon_{1}&=&\Upsilon_0\cup P(\Upsilon_0)\cup\Upsilon_0P(\Upsilon_0),\\
\vdots\ \ & &\ \ \ \ \vdots\\
\Upsilon_{n}&=&\Upsilon_0\cup P(\Upsilon_{n-1}) \cup \Upsilon_0P(\Upsilon_{n-1}), \\
\vdots\ \ & &\ \ \ \ \vdots
\end{eqnarray*}
and
$$
\Upsilon(D^{\omega}(X))=\bigcup_{n\geq 0 }\Upsilon_n.
$$

\begin{theorem}\label{t6.2}
$Irr(S)= \Upsilon(D^{\omega}(X))$  is a basis of  $K[ X;\Omega|S]$.
\end{theorem}
{\bf Proof:} By  Theorem \ref{3.4} and Theorem \ref{t6.1}, we can
obtain the result easily. \hfill $\blacksquare$

\begin{theorem}
$K[ X;\Omega|S]$ is a free commutative $\lambda$-differential
Rota-Baxter algebra on  set $X$ with a basis
$\Upsilon(D^{\omega}(X))$.
\end{theorem}
{\bf Proof:}  By Theorem \ref{t6.2}. \  \

\end{document}